\documentclass[letterpaper, 10 pt, conference]{ieeeconf}  

\IEEEoverridecommandlockouts                              

\usepackage{color}
\usepackage{hyperref}
\usepackage{amsmath} 
\usepackage{bbm}
\usepackage{amssymb}
\usepackage{lipsum}
\usepackage{breqn}
\usepackage{cancel}
\usepackage{float}
\hypersetup{
    colorlinks=true,
    linkcolor=blue,
    filecolor=magenta,      
    urlcolor=cyan,
    }
\newtheorem{defn}{Definition}
\newtheorem{thm}{Theorem}
\newtheorem{lemma}{Lemma}

\newtheorem{prop}{Proposition}

\usepackage{algorithm}
\usepackage{bm}
\usepackage{algpseudocode}

\usepackage{caption}

\newcommand{\E}{\mathbb{E}}
\newcommand{\R}{\mathbb{R}}
\newcommand{\Z}{\mathbb{Z}}
\newcommand{\1}{\mathbf{1}}

\newcommand{\Wv} {{\bf W}}
\newcommand{\BiObj}{B} 
\newcommand{\BiCon}{C}

\def\one{\mbox{1\hspace{-4.25pt}\fontsize{12}{14.4}\selectfont\textrm{1}}} 
\DeclareMathOperator*{\argmin}{\arg \min} 

\newcommand{\TotCostT}[5]{v^{#1}_{{#2},{#3}}({#4},{#5})} 

\newcommand{\ignore}[1]{}  

\title{\LARGE \bf
Finite-Horizon Constrained MDPs \\ With Both Additive And Multiplicative Utilities
}

\author{Uday Kumar M$^{1}$, Sanjay P Bhat$^{1^*}$, Veeraruna Kavitha $^{2}$ and Nandyala Hemachandra $^{2^*}$
\thanks{$^{1,1^*}$ Uday Kumar M and Sanjay P Bhat are with TCS Research, Hyderabad, India
        {\tt\small \{udaykumar1.m,sanjay.bhat\}@tcs.com}}%
\thanks{$^{2,2^*}$ V. Kavitha and N. Hemachandra are faculty members with Dept of Industrial Engineering and Operations Research, Indian Institute of Technology Bombay, Mumbai, India
        {\tt\small \{vkavitha,nh\}@iitb.ac.in}}
}

\begin{document}

\maketitle
\thispagestyle{empty}
\pagestyle{empty}

\begin{abstract}
This paper considers the problem of finding a solution to the finite horizon constrained Markov decision processes (CMDP) where the objective as well as constraints are sum of additive and multiplicative utilities. Towards solving this, we construct another CMDP, with only additive utilities under a restricted set of policies, whose optimal value is equal to that of the original CMDP. 
Furthermore, we provide a finite dimensional bilinear program (BLP) whose value  equals the CMDP value and whose solution provides the optimal policy. We also suggest an algorithm to solve this BLP.
\end{abstract}

\section{INTRODUCTION}

\ignore{\color{red}Markov decision processes (MDPs) form a popular and widely-used framework for dynamic optimization in a stochastic control setting. The literature on the theory and applications of standard MDPs involving discounted, total or averaged additive costs is extensive and goes back several decades \cite{PUT,Arnab_Jeffrey}. On the other hand, risk-sensitive MDPs (RSMDPs) involving expected exponential utility of sum of stage-wise costs (\cite{UdaySanjayKavitaNH,SJaq1973,Kavitha_NH_Atul}) are relatively less studied. One possible reason for this could be that optimal policies in the case of RSMDPs tend to be non-stationary, and efficient solution methods are not available in the case where multiple constraints are present. RSMDPs form a special case of MDPs where rewards and costs are multiplicative  which involve expected products of stage-wise costs. Such MDPs with multiplicative reward  have been studied \cite{Karel} to an even lesser extent than RSMDPs.}

In this paper, we consider a class of constrained Markov Decision Processes   (CMDPs) which reduce to finite-horizon standard MDPs, risk-sensitive MDPs as well as multiplicative MDPs as special cases. More specifically, we consider a finite-state CMDP whose objective as well as constraints involve a linear combination of additive and multiplicative costs.   Multiplicative costs such as risk-sensitive costs  provide  one way of achieving  robust optimization by optimizing a weighted sum of moments (\cite{UdaySanjayKavitaNH, SJaq1973, Kavitha_NH_Atul}). On the other hand 
additive  components model a variety of costs, like budget-related,  and we are interested in optimizing the expected value of a weighted combination of the two types of costs such as $    
\sum_{t=0}^T c_t(X_t,A_t) + \exp\left(\gamma\sum_{t=0}^T r_t(X_t,A_t)\right).
$
\ignore{
For example, in a  queuing system one might want to optimize a multiplicative robust cost related to waiting time or loss fractions (as in \cite{Kavitha_NH_Atul}) along with the {\color{red}{cumulative ??}} average server capacity utilized over the given time horizon.} 

Unconstrained MDPs with purely multiplicative cost components can be solved using dynamic programming \cite{Kavitha_NH_Atul}. Reference \cite{Kavitha_NH_Atul}  gave another solution that is LP-based technique for a finite-horizon MDP having a purely multiplicative (specifically, risk-sensitive) objective cost and at most one purely additive constraint cost. The solution technique of \cite{Kavitha_NH_Atul} is based on augmenting the state space with a discrete variable that keeps track of the running multiplicative cost along the path. As a result, the LP provided in \cite{Kavitha_NH_Atul} grows exponentially in size with the time horizon. Moreover, the technique used in \cite{Kavitha_NH_Atul} does not extend to the case where one or more constraints involve multiplicative costs. 
A recent algorithm in \cite{Varthika_Kavitha}  solves constrained risk sensitive MDPs, however they consider only  risk-sensitive components. 
In this paper, we provide a Bilinear programming (BLP) based technique for solving 
the CMDP problems with both multiplicative and additive costs. 
\ignore{
the considerably more difficult and hitherto unattempted problem involving an CMDP with multiple constraints, where the objective and constraints involve a linear combination of additive and multiplicative components. } 

Our BLP based solution technique involves augmenting the state with binary variables, one for every multiplicative component appearing in either the objective or the constraints. We construct a CMDP defined on the augmented state space such that the additive and multiplicative cost components of the original CMDP are absorbed into the additive stage-wise costs and the controlled transition functions, respectively,  of the augmented MDP. 

We further suggest an  iterative algorithm to solve the BLP  (inspired by  \cite{Varthika_Kavitha}). We conclude this section with an application of the CMDP considered in this paper.

\subsection*{Epidemics and Delay Tolerant Networks}
\ignore{
There are many applications,    where  one encounters combined cost like in \eqref{eq_obj_comb_const}-\eqref{eq_constr_comb} example, risk-sensitive cost is one way of considering  robust optimization, where one attempts to optimize weighted sum of moments (\cite{UdaySanjayKavitaNH, SJaq1973, Kavitha_NH_Atul}).
But there can be other cost components, like budget related,  and agents might be interested in optimizing a weighted combination of such costs
$$
\sum\limits_{t=0}^T c_t(X_t,A_t) + e^{\gamma\sum\limits_{t=0}^T r_t(X_t,A_t)}
$$
For example in a  queuing system one might want to optimize a robust cost  related to waiting time/loss fractions (as in \cite{Kavitha_NH_Atul}) along with the average server capacity utilized over the given time frame. 
One can encounter many such applications. 

One can also encounter other set of applications, where the objective function itself is risk sensitive. We would explain one such scenario in the following.}

Consider an area with $N$ users
where an epidemic is fast spreading. The government has to devise a lock-down strategy to efficiently combat the disease, keeping in view of lock-down costs (in terms of economical losses). The   time-frame is divided into $T$-time slots and the policy is to design the level of  lock-down to be imposed for each time slot, based on the system state at the beginning of the same slot. 

Any infected user can infect a normal/susceptible user
when the former comes in contact with the latter (see \cite{Varthika_Khushboo_Kavitha} for similar details). The successive contacts between any two users are modelled by a Poisson processes (as in \cite{DTNs,Varthika_Khushboo_Kavitha}). The  rate of this contact process $\lambda$ is proportional to the imposed level of lock-down.
Let $\Lambda_t$ represent contact rate chosen in slot $t$ and let  $g(\Lambda_t)$ be the economic cost corresponding to the related lockdown level.      

Let  $X_t$  represent system state,   the number of infected individuals at the beginning of   time slot $t$;  a person infected in a time slot can infect from    next time slot onwards. 
Then 
$$
X_{t+1} = {\cal B} (N-X_t,  q_t) - {\cal B} (X_t, r), 
$$where ${\cal B} (\cdot, \cdot)$ is a binomial random variable,  $r$ is the probability that an infected user gets recovered in a time slot and  $q_t = 1-\exp (-\Lambda_t X_t)$ is the probability that a susceptible gets infected in slot $t$. 
We are interested in the probability $P_S (\pi)$ that a given typical agent survives without infection in the given time-frame for a given lock-down policy $\pi$. 
By conditioning on state trajectory $\{X_t\}$ (details as in \cite{DTNs})
$
P_S(\pi)  = \E^\pi \left [ e^{-\sum_{t=1}^T \Lambda_t X_t } \right ] 
$. Thus in all we would optimize a combined cost that also considers the economic losses due to lock-down strategy as   below  ($\alpha$ -- trade-off factor):
$$
\min_\pi \E^\pi \left [\sum_t g(\Lambda_t) + \alpha e^{-\sum_{t=1}^T \Lambda_t X_t } \right ].
$$ 
\ignore{Further there is  a constraint on the expected number of infected individuals; by conditioning on $X_t$, this constraint is given  as below for some $b < \infty$:
$$\E^\pi\left [\sum_t {\cal B} (N-X_t,  q_t) \right ] = E\left [ \sum_t (N-X_t) q_t \right ] \le b.$$}

A similar problem  (excluding recoveries) arises  in delay tolerant networks (\cite{DTNs}), when one considers a combined cost   related to  delivery failure probability and power-budget. One can further have a bound on the expected number of copies. 

The aim of this paper is to  consider a general class of problems that are of the nature as in the above example.

\ignore{are moving and consider propagation of infection when a healthy (or susceptible) comes in contact with an infected person. One might be interested in the probability $P_S$ that a given typical agent survives without infection in the given time-frame, given the contact process. One might be able to control the contact process. This for example could represent an epidemic scenario in which case one wants maximize the above infection-survival probability $P_S$, or could represent DTN scenario in which one might want to maximize the infection probability $(1-P_S)$, which in this case represents the probability of successful message  transmission to the given destination (details in \cite{DTN-Kavitha}). }

\begin{figure*}[ht]
Let $\Wv_1 := \{ w_t (x, \1, a) ; \mbox{ for all } x, a, t<T\}$ and 
$\Wv_2 := \{ w_t (x, z, a) ; \mbox{ for all } x, a, t<T, z \ne \1 \}, \Wv = (\Wv_1, \Wv_2).$
\hrule 

\begin{dmath}
\mbox{\bf BLP ($ \Wv_1, \Wv_2$): } 
{\min\limits_{\{\Wv = (\Wv_1, \Wv_2) \} } \sum\limits_{t=0}^{T-1}\sum\limits_{(x,z)\in\mathcal{\bar{X}}}\sum\limits_{a\in\mathcal{A}}\bar{r}_{t,0}((x,z),a) 
 w_{t}((x,z),a)
 }
 +  \sum\limits_{(x,z)\in\mathcal{\bar{X}}} \bar{r}_{T,0}(x,z)w_{T}(x,z)\\
\mbox{ sub to.} {\sum\limits_{a\in\mathcal{A}} w_{0}((x,z),a) = \one_{\{(s,\1)\}}(x,z),  \   \mbox{ all the constraints  are for all } (x,z)\in\mathcal{\bar{X}},}  \mbox{ and } t \in \{1 \cdots T\} \\ 
\sum\limits_{a\in\mathcal{A}} w_{t}((x,z),a)  - \sum\limits_{(x^\prime,z^\prime)\in\mathcal{\bar{X}}}\sum\limits_{a\in\mathcal{A}}\bar{Q}_t\bigg((x,z)\bigg|(x^\prime,z^\prime),a\bigg){w_{t-1}((x^\prime, z^\prime),a)=0},\\ 
{ w_{T}(x,z)-\sum\limits_{(x^\prime,z^\prime)\in\mathcal{\bar{X}}}\sum\limits_{a\in\mathcal{A}}\bar{Q}_{T}\bigg((x,z)\bigg|(x^\prime,z^\prime),a\bigg)}
 {w_{T-1}((x^\prime,z^\prime),a)= 0,} 
 \\
{\sum\limits_{t=0}^{T}\sum\limits_{(x,z)\in\mathcal{\bar{X}}}\sum\limits_{a\in\mathcal{A}}\bar{r}_{t,i}((x,z),a) w_{t}((x,z),a)  
 +  \sum\limits_{(x,z)\in\mathcal{\bar{X}}}  \bar{r}_{T,i}(x,z)w_{T}(x,z)\leq b_i  \condition{for } i=1,\ldots, K,}\\
{w_t(x, z, a) \sum_{a'} w_t (x, \1, a') = w_t (x, \1, a) \sum_{a'} w_t (x, z, a'),
\  \ 
w_{t}((x,z),a)\geq 0, w_{T}(x,z)\geq 0,   \forall a\in\mathcal{A}.}\\
\label{eq_LP}
\end{dmath}

\vspace{3mm}
\hrule
\end{figure*}

\section{Model Formulation and Problem Statement}
We consider a finite-horizon Markov decision process (MDP) $M:=(\mathcal{X},\mathcal{A},Q,\{0,1,\ldots,T\})$, where $\mathcal{X}$ is the finite state space, $\mathcal{A}$ is the finite action space, 
$T>0$ is the terminal time, $Q$ is the transition function (or transition law) and $\{0,1,\cdots,T\}$ is the set of discrete decision epochs.
Here
$Q(x^\prime|x,a)$ is the probability that the system transitions to state $x^\prime$ when action $a$ is taken at state $x$. 

A Markovian randomized (MR) decision rule is a map $d:\mathcal{X}\to \mathcal{P(A)}$, where $\mathcal{P(A)}$ is the set of probability distributions on the action space $\mathcal{A}$. We denote by $d(a|x)$ the probability of choosing action $a$ in the state  $x$ under the MR decision rule $d$. An MR policy is a sequence of MR decision rules indexed by the decision epochs. We denote the set of MR policies by $\Pi_{\rm MR}$. For every initial state $s$, a policy $\pi:=\{d_t\}_{t=0}^{T-1
}\in\Pi_{\rm MR}$  induces a probability measure $P^{\pi}_s$ on the space of state-action trajectories (see \cite[Ch2.]{PUT}).
\ignore{$\Omega$  (Chapter 2  of \cite{PUT} and     \cite[Prop C.10, Rem C.11 and Ch 2]{LermaLasserre}),
\begin{eqnarray}
P^\pi_s(X_0=x) = \one_{\{x=s\}}\label{eq_init_distX}
\end{eqnarray}
\begin{eqnarray}
&P^\pi_s(A_t=a|X_t=x,X_{t-1}=x_{t-1},A_{t-1}=a_{t-1}, \nonumber \\
&\ldots, X_{0}=x_{0},A_{0}=a_{0})= d_t(a|x)\label{eq_prob_actX}
\end{eqnarray}
\begin{eqnarray}
&P^\pi_s\bigg(X_t=x^\prime\bigg|X_{t-1}=x,A_{t-1}=a,\ldots,X_{0}=x_{0},A_{0}=a_{0}\bigg) \nonumber \\
&= Q\bigg(x^\prime\bigg|x,a\bigg)\label{eq_policy_transX}
\end{eqnarray}
The construction of such measure is by defining the probability of the trajectory by
\begin{eqnarray}
&P^\pi_s(x_0,a_0,x_1,a_1,\ldots,x_{T-1},a_{T-1},x_T)\nonumber\\
&=\one_{\{x_0=s\}}d_0(a_0|x_0)Q(x_1|x_0,a_0)\nonumber\\
&\ldots d_{T-1}(a_{T-1}|x_{T-1})Q(x_T|x_{T-1},a_{T-1}).\label{eq_Prob_Traj_pi}
\end{eqnarray}
}
 Let $\E^\pi_s[\cdot]$ denote the corresponding expectation. 

This paper considers an MDP involving additive as well as multiplicative stage-wise components; we represent these  by $r_{t,i}:\mathcal{X}\times \mathcal{A}\times \mathcal{X}\rightarrow \mathbb{R}$ and $f_{t,i}:\mathcal{X}\times \mathcal{A}\times \mathcal{X}\rightarrow \R$, respectively,  where $0\leq i\leq K$ and $t\in\{0,1,\ldots,T-1\}$. More precisely, we consider the following optimization problem: given $s\in\mathcal{X}$ and $b_i\in\R$ for each $0\leq i\leq K$, find
\begin{equation}
    \inf_{\pi\in\Pi_{\rm MR}}w^\pi_0(s), \ \mbox{subject to, } w^\pi_i(s) \leq b_i, \ \forall \ 1\leq i\leq K, \label{eq_CMDP_Prob} \tag{$P$}
\end{equation}
where, for any policy $\pi\in\Pi_{\rm MR}$, for each $0\leq i \leq K$ and $\alpha_i\in\R$, we define

\vspace{-3 mm}
{\small
\begin{eqnarray*}
 w^\pi_i(s) \stackrel{\triangle}{=}\E^\pi_s\bigg[\sum\limits_{t=0}^{T-1} r_{t,i}(X_t,A_t,X_{t+1})+\alpha_i \prod\limits_{t=0}^{T-1}f_{t,i}(X_t,A_t,X_{t+1})\bigg].
\end{eqnarray*}}


\ignore{
We introduce stage-wise cost functions  $r_{t,i},f_{t,i}:\mathcal{X}\times\mathcal{A}\times \mathcal{X}\to\R$ for each $t\in\{0,1,\ldots,T-1\}$ and $i\in\{0,1,\ldots,K\}$. Let $\alpha_i,b_i\in\R$ be fixed for $i\in\{0,\ldots,K\}$. For any policy $\pi$, $i\in\{0,\ldots,K\}$ and initial state $s$, define,
\begin{equation}\begin{split}
w^\pi_i(s) &:=\E^\pi_s\bigg[\sum\limits_{t=0}^{T-1} r_{t,i}(X_t,A_t,X_{t+1}) \nonumber\\ 
&\quad+\alpha_i \prod\limits_{t=0}^{T-1}f_{t,i}(X_t,A_t,X_{t+1})\bigg].\label{eq_wpi_def}
\end{split}
\end{equation}
 Let $\mathcal{F}:= \{\pi\in\Pi_{\rm MR}: w^\pi_i(s)\leq b_i \mbox{, } \forall 1\leq i\leq K\}$ the set of policies that is of interest to us. This paper attempts to solve the following constrained optimization problem:
\begin{equation}
\min\limits_{\pi\in\mathcal{F}}w^\pi_0(s).    \label{eq_CMDP_Prob}
\end{equation}}

We assume that, for each $t$ and $i$ 
the multiplicative cost component $f_{t,i}$ has the same sign everywhere on $\mathcal{X}\times\mathcal{A}\times \mathcal{X}$. 
By appropriately scaling and then absorbing the scaling factor and common sign into the coefficient $\alpha_{i}$, we  assume  
without  loss of generality that $0\leq f_{t,i}\leq 1$ for all $t$ and $i$. 

\newcommand{\eop}{\hfill $\blacksquare$}
\subsection{Comparison and Special Cases}

For finite horizon problem, one may argue that the problem can be converted to standard linear cost MDP problem,  by augmenting the accumulated multiplicative cost  as  an additional state, one for each multiplicative component, e.g., $Y_{\tau,i} := \prod_{t=0}^{\tau-1} f_{t,i} (X_t, A_t, X_{t+1}) $.
The complexity of  problem obtained after such an augmentation is extremely high; problem \eqref{eq_CMDP_Prob_Augm} (our solution defined in later sections) also augments one state per multiplicative component, but the augmented state component is always binary valued, as opposed to the geometric increase  with $t$ in size of the state space that accommodates $\{Y_{t,i}\}$. 
%
%

The problem defined in \eqref{eq_CMDP_Prob} is fairly general, and covers many  special cases listed below.  

\noindent$\bullet$ With $\alpha_i = 0$ for all $i$, we have the well known standard discounted cost MDP with constraints (\cite{Eitan_CMDP}). 


\noindent$\bullet$ By setting $r_{t,i}=0$ and $f_{t,i}(x,a,x^\prime)=\exp(\beta^{t}c_i(x,a))$   the problem \eqref{eq_CMDP_Prob} is the well known    {\it risk-sensitive MDP} (\cite{UdaySanjayKavitaNH, Kavitha_NH_Atul}), solutions are seldom known for such problems (except for \cite{Varthika_Kavitha}) and we are providing an easily implementable solution.  


\noindent$\bullet$  The state space may contain some undesirable  {\it ``error" states}  denoted by $\mathcal{E}$ (\cite{PeterFritz}). The problem of optimizing the cost with constraint that the probability of entering the error state is bounded by a threshold can be handled   by setting $f_{t,1}(x,a,x^\prime)=\one_{\mathcal{E}}(x,a)$ and $r_{t,i} \equiv 0$  in \eqref{eq_CMDP_Prob}.

\noindent$\bullet$ In addition, one can also model applications like the one described in Introduction with combined costs.  

\section{Equivalent CMDP And Main results}
 In this section, we provide the two main results of this work. The first result gives the equivalence between  the original CMDP \eqref{eq_CMDP_Prob} and a newly constructed CMDP involving only additive costs. The second result provides  an equivalence to a finite dimensional optimization problem with linear objective function and with linear and bilinear constraints.  The proofs of the results are  in the appendix.

We augment the state space with one binary variable per multiplicative cost, which starts with all 1s and absorbed to 0. The core idea is to capture the multiplicative cost via the expectation of terminal value of the binary variable. In the next subsection, we give precise mathematical  details of this construction.

\subsection{Equivalent CMDP}
We construct a new CMDP involving only additive cost and constraint components. The new CMDP is defined by augmenting the original state with one binary variable for each non-zero $\alpha_{i}$ in \eqref{eq_CMDP_Prob}. To simplify the exposition, suppose $\alpha_{i}\neq 0$ for each $i$. Define the augmented state space by $\mathcal{\bar{X}}:=\mathcal{X}\times \mathcal{Z}$, where $\mathcal{Z}= \{0,1\}^{K+1}$. For each $(x,z),(x^\prime,z^\prime)\in\mathcal{\bar{X}}$,  $a\in\mathcal{A}$ and $t\in\{0,\ldots,T-1\}$ the transition function of the new CMDP at   epoch $t$ is   defined by  
\begin{equation}\begin{split}
&\bar{Q}_t\bigg((x^\prime,z^\prime)\bigg|(x,z),a\bigg)\\
&=P((X_t,Z_t)=(x^\prime,z^\prime)|(X_{t-1},Z_{t-1})=(x,z),A_{t-1}=a))\\
& :=\begin{cases}
0, \mbox{ if } z_i=1-z_i^\prime=0 \mbox{ for at least one }  i,\\
Q(x^\prime |x,a)\prod\limits_{i=0}^{K} \rho_{t-1,i}(x,x^\prime,a,z,z^\prime),\mbox{ otherwise},
\end{cases}\nonumber
\end{split}
\end{equation}
%
\begin{eqnarray}
\mbox{where } \rho_{t,i}(x,x^\prime,a,z,z^\prime) &:=&  (f_{t,i}(x,a,x^\prime))^{z_iz_i^\prime}\nonumber\\ 
&&\times (1-f_{t,i}(x,a,x^\prime))^{z_i(1-z_i^\prime)},\nonumber
\end{eqnarray}
for each $t\in\{0,\ldots,T-1\}$ and $0\leq i\leq K$.

It is easy to check that the transition law $\bar{Q}_t(\cdot|(x,z),a)$ is indeed a probability distribution on the set $\mathcal{\bar{X}}$ for all $(x,z)\in\bar{X}$ and $a\in\mathcal{A}$. This allows us to define a MCM, $\bar{M}:=(\mathcal{\bar{X}},\mathcal{A},\bar{Q}_t,\{0,1,\ldots,T\})$ to represent the state evolution of the  process $\{\bar{X}_t\}_{t=0}^{T}$ where $\bar{X}_t:=(X_t,\{Z_{t,i}\}_{i=0}^{K})$. Here, for every $i$, $Z_{t,i}\in\{0,1\}$. The central idea of this construction is two fold: (a) the transitions of the original Markov chain $\{X_t\}_{t=0}^{T}$  are not affected by the process $\{Z_t\}_{t=0}^{T}$ and (b) as we shall soon see, the expected value of  the binary variable  $Z_{T,i}$ equals the expected value of the multiplicative cost $\prod_{t}f_{t,i}(X_{t},A_{t},X_{t+1})$.

Next, we define $\bar{\Pi}$ to be the set of MR policies w.r.t $\bar{M}$ which are indifferent to values of the augmented state. More precisely, $\bar{\Pi}$ is the set of policies  $\bar{\pi}=\{\bar{d}_t\}_{t=0}^{T-1}$ such that $\bar{d}_t(a|(x,z))=\bar{d}_t(a|(x,\1))$ for all $t\in\{0,\ldots,T-1\}$ 
where $\1:=(1,\ldots,1)\in\R^{K+1}$. It is easy to observe that $\bar{\Pi}$ is in a one-to-one correspondence with $\Pi_{\rm MR}$, the space of MR policies w.r.t the original model $M$. We set the initial state in $\bar{M}$ to be $(s,\1)$, where $s$ is initial state in $M$.


Finally, we define the stage-wise costs in $\bar{M}$ by
\begin{eqnarray}
\bar{r}_{t,i}((x,z),a)&:=&\sum\limits_{x^\prime\in\mathcal{X}} r_{t,i}(x,a,x^\prime)Q(x^\prime|x,a),\nonumber\\
 \bar{r}_{T,i}(x,z)&:=&\alpha_i z_i \  \mbox{ for all } i, t<T.
\end{eqnarray}
Finally, we define the following optimization problem the newly constructed augmented MDP,
\begin{equation}\label{eq_CMDP_Prob_Augm}\tag{$\bar{P}$} 
        \begin{split}&\min\limits_{\eta\in\bar{\Pi}} \TotCostT{\eta}{T}{0}{s}{\1} \\
        \mbox{ subject to, } &\TotCostT{\eta}{T}{i}{s}{\1}\leq b_i,\mbox{ }\forall ~ 1\leq i\leq K, \mbox{ where},
        \end{split}
\end{equation}
\vspace{-0.3 cm}
{\small
\begin{eqnarray*}
\TotCostT{\eta}{T}{i}{x}{z} := \E^{\eta}_{(x,z)}\bigg[\sum\limits_{t=0}^{T-1}\bar{r}_{t,i}\bigg((X_t,Z_t),A_t\bigg)  + \bar{r}_{T,i}(X_T,Z_T)\bigg]. 
\end{eqnarray*}}
Observe, that the objective and constraints of the problem \eqref{eq_CMDP_Prob_Augm} are both linear/additive only. We define the projection map $\Gamma:\bar{\Pi}\mapsto \Pi_{\rm MR}$ as follows; for any $\eta=\{\bar{d}_t\}_{t=0}^{T-1}\in\bar{\Pi}$, define $\Gamma(\eta) = \{d_t\}_{t=0}^{T-1}$ where $d_t(a|x)=\bar{d}_t(a|x,\1)$ for all $t$, $a$ and $x$.
\subsection{Main Results}
We now give the first main result of this paper: the equivalence of \eqref{eq_CMDP_Prob} and \eqref{eq_CMDP_Prob_Augm}. 
\begin{thm}\label{thm_main_thm1}
The optimal value of the two problems \eqref{eq_CMDP_Prob} and \eqref{eq_CMDP_Prob_Augm} are equal.
Furthermore, $\eta^* $ is an optimizer  of \eqref{eq_CMDP_Prob_Augm} if and only if $\Gamma(\eta^*)$ is optimizer of \eqref{eq_CMDP_Prob}.      $\blacksquare$

\end{thm}

 Theorem \ref{thm_main_thm1} provides an alternative CMDP \eqref{eq_CMDP_Prob_Augm} where the objective as well as constraints are additive (or linear) as in a standard MDP, but the policy space is restricted to the policies which are indifferent to the augmented state component. 
 There many solution techniques to solve  standard MDP  (\cite{Eitan_CMDP}), however the restricted policy space requires special attention.

%
The second main result of this paper is  the Bi-linear programming  (BLP), given by \eqref{eq_LP} provided at the top of the page, which   solves the new problem \eqref{eq_CMDP_Prob_Augm}.  The bilinear constraint in BLP \eqref{eq_LP} ensures that the solution is indifferent to augmented state $z$.  
\begin{thm}{\bf [Solution using BLP]}\label{thm_main_thm2}
 The value of  BLP \eqref{eq_LP} is the   value of the problems \eqref{eq_CMDP_Prob} and \eqref{eq_CMDP_Prob_Augm}. Let $w_{t}^*((x,z),a)$ be the solution of the BLP. Then the optimal policies $\bar{\pi}^*:=\{\bar{d}^*_t\}_t$ for the augmented problem \eqref{eq_CMDP_Prob_Augm} and $\pi^*:=\{d^*_t\}_t$ for the original \eqref{eq_CMDP_Prob} are given by
 \begin{eqnarray}
  \bar{d}_t^*(a|(x,z)) &=& d_t^*(a|x) = \frac{w_{t}^*((x,\1),a)}{\sum\limits_{a^\prime\in\mathcal{A}}w_{t}^*((x,\1),a^\prime)},\nonumber\\
  &&\mbox{ for }x\in\mathcal{X}, \mbox{ } a\in\mathcal{A}, z\in\mathcal{Z}, t<T.\label{eq_sol_orig_prb} 
 \end{eqnarray}
\end{thm}

\renewcommand{\S}{\mathcal{S}}
\newcommand{\A}{\mathcal{A}}
\newcommand{\bbtheta}{{\bm \Theta}}
\newcommand{\btheta}{{\bm \theta}}
\newcommand{\bQ}{{\bf Q}}
\newcommand{\bbQ}{{\bm {\mathbbm Q}}}
\newcommand{\B}{\mathcal{B}}
\newcommand{\F}{\mathcal{F}}
\newcommand{\bB}{\bm{\mathcal{B}}}
\newcommand{\bPi}{\Wv}
\newcommand{\bF}{\bm{\mathcal{F}}}
\newcommand{\M}{\mathcal{W}}
\newcommand{\Mset}{\mathbbm{W}}
\newcommand{\cP}{\mathbbm{P}}
\newcommand{\tpi}{{\tilde{\bPi}}}
\newcommand{\mdp}{\mbox{combined-CMDP}}
\newcommand{\rpi}{{\Psi}}
\newcommand{\lp}[1]{\mbox{Bi-LP}(#1)}

\newcommand{\LP}{{\cal L}}
 
\section{Algorithm}

 By Theorem \ref{thm_main_thm2}, 
   constrained global optimization problem of BLP   \eqref{eq_LP} is equivalent to the $\mdp$ problem given in  \eqref{eq_CMDP_Prob_Augm}. We now suggest an algorithm to solve the BLP which can rewritten as  
\begin{equation}
\inf_{\Wv} \BiObj (\Wv) \mbox{ s.t. } \BiCon_1 (\Wv) \le {\bf b} \mbox{ and } \BiCon_2 (\Wv) = 0,\tag{$BLP$}\label{eq_BLP_rewritten}    
\end{equation}
for appropriate functions $(\BiObj, \BiCon_1, \BiCon_2)$ and  ${\bf b}$, for example:
\begin{eqnarray}
\label{Eq_BiObj}
\BiObj (\Wv) &=& \sum\limits_{t=0}^{T-1}\sum\limits_{(x,z)\in\mathcal{\bar{X}}}\sum\limits_{a\in\mathcal{A}}\bar{r}_{t,0}((x,z),a) 
 w_{t}((x,z),a) \nonumber \\
 &&+  \sum\limits_{(x,z)\in\mathcal{\bar{X}}} \bar{r}_{T,0}(x,z)w_{T}(x,z).
 \end{eqnarray}
This BLP has $\Wv = (\Wv_1, \Wv_2)$ as variables, and when one of them is fixed then it  is clearly an LP in the other variable. 
Let $\LP_1 (\Wv_2)$ and $\LP_2 (\Wv_1)$ represent the solution set of the respective LPs when $\Wv_2$ and $\Wv_1$ are fixed. 
One can define a relevant 
  fixed point equation using  the above two LP solutions, to be precise we are interested in the fixed points $\bPi \in \M (\bPi)$, where  the solution set $\M$ is defined as,
\begin{eqnarray}
\M (\bPi):=  \{ \tpi : \tpi_1 \in \LP_1 (\bPi_2), \tpi_2 \in \LP_2 (\bPi_1)   \}.\label{eq_sol_set_LP12}    
\end{eqnarray}

Define $\Mset = \{\bPi : \bPi \in \M (\bPi)\}$, the set of all fixed points.  It is not difficult to verify that  the solution of the following global fixed point problem \eqref{eq_FP_optim_problem} provides the solution to BLP \eqref{eq_LP}:
\begin{equation}
\inf_{\Wv \in \Mset}  \BiObj (\Wv). \tag{$GF$} \label{eq_FP_optim_problem}  
\end{equation}
Thus any fixed point iterative algorithm is useful in solving the BLP, however one needs to ensure that it is converging towards the best fixed point (as define above).  The global optimization problem defined in \eqref{eq_BLP_rewritten} can be useful in this context. This solution approach provided below, exactly parallels that in the recently provided algorithm \cite{Varthika_Kavitha}; the authors  in \cite{Varthika_Kavitha} consider constrained risk-sensitive MDP problem and also have a fixed point equation  and global optimization problem as in \eqref{eq_BLP_rewritten} and \eqref{eq_FP_optim_problem}. We believe the justification of the algorithm can parallel that provided in  in \cite{Varthika_Kavitha}, however skip those details due to lack of space.

 This global optimization problem can for example be solved   using   random restarts \cite{Pepelyshev_etal}.  It  has two types of update steps (at any iterative step $k$): i) \underline{random search step} -- a random new point is chosen  from   the feasible region 
   with   probability $p_k$, and, ii) \underline{local improvement step} as in \eqref{eqn_local_update} is chosen otherwise. 
 The probability $p_k$ diminishes with $k$.

The aim in the local improvement step is to converge to a fixed point in $\Mset$ as given by the following:
\begin{eqnarray}\label{eqn_local_update}
    \bPi_{k+1}& = &  \bPi_k+ \epsilon_k (\rpi_{k}(\bPi_k) -\bPi_k),  \\
    \mbox{ where, }\rpi_k(\bPi_k) &\in& \M(\bPi_k) \mbox{ and, }\nonumber 
    \epsilon_k  = 1/k. \nonumber
\end{eqnarray}In the above equation, $\rpi_{k}(\bPi_k)$ is chosen randomly from $\M(\bPi_k)$, the solution set  \eqref{eq_sol_set_LP12} derived using the two LPs of $\lp{\bPi_k}$;  basically $\bPi_k = (\bPi_{k,1}, \bPi_{k,2})$ and, $$\rpi_{k,1}(\bPi_{k}) \in LP_1(\bPi_{k,2}) \mbox{ and  } \rpi_{k,2}(\bPi_{k}) \in LP_2(\bPi_{k,1}).$$

\begin{algorithm}[h]
\caption{Global $\mdp$ algorithm (GRC)}\label{alg_random_search}
Initialize $\bPi_0$ randomly,  set $\BiObj^*= -\infty$, $\hat{\bPi}^*=\bPi_0$; choose a constant $w$; ${\cal U}$ set of all possible policies.

\textbf{For} $k=1,2,\dots$
\begin{algorithmic}
\State \vspace{-5mm}\begin{eqnarray*}
\bPi_k \hspace{-1mm} \gets \hspace{-1mm} \left\{ \hspace{-1mm}
\begin{array}{ll}
     \mbox{random policy chosen from } {\cal U}&  \mbox{w.p. } p_k = \frac{w}{k} \\
     \mbox{Local improvement }  (\bPi_{k-1}) \mbox{ of } \eqref{eqn_local_update}&  \mbox{w.p. } 1-p_k \\ 
\end{array}\right.
\end{eqnarray*}
\State Calculate $\BiCon_1(\bPi_k)$ and $\BiCon_2(\bPi_k)$  

\If{$\BiCon_2(\bPi_k) =0$ and $ \BiCon_1(\bPi_k) \le  {\bf b}$} 
\State Calculate $\BiObj(\bPi_k)$ using \eqref{Eq_BiObj}
\vspace{-4mm}
\State \If{$\BiObj(\bPi_k) \le \BiObj^*$}
\State $\hat{\bPi}^*\gets \bPi_k$
\State $\BiObj^* \gets \BiObj(\bPi_k)$
\EndIf
\Else 
\State Choose random policy  from ${\cal U}$   \   (random restart again)
\EndIf
\end{algorithmic}
\end{algorithm}

\subsection*{Complexity comparison}
The paper \cite{Kavitha_NH_Atul} gives the LP based solution to risk-sensitive objective and a single linear additive constraint. 
The table \ref{tab_constr_compar} compares the  complexity of one constraint CMDP solved using LP given in \cite{Kavitha_NH_Atul} and our BLP. We observe that our BLP is linear in $T$ while the LP of \cite{Kavitha_NH_Atul} is exponential in $T$ that makes the latter more complex to solve the problem.
\begin{table}[h!]  
\begin{center}  
\caption{No. of Variables \& Constraints }
\label{tab_constr_compar}  
\begin{tabular}{|c|c|c|c|}  
\hline
\textbf{Reference} & \textbf{No of Variables}  &\textbf{No of Constraints}  \\  
\hline
LP as in \cite{Kavitha_NH_Atul} & $\frac{mn ((mn)^T-1)}{(mn-1)}$ & $m+ \frac{m((mn)^T-1)}{(mn-1)} + 1 $\\\hline
BLP &$4m(Tn+1)$ &$4m(2+2T+Tn)+1$  \\\hline
\end{tabular}
\end{center}
\end{table}

\ignore{
\begin{figure*}
    \centering
\includegraphics[width=0.9\textwidth, height=0.8\textwidth]{Discussion_with_Prof_Kavaith.jpg}
    \caption{Caption}
    \label{fig:my_label}
\end{figure*}
}

\section{Proof of theorems}

The domains of the problems \eqref{eq_CMDP_Prob} and \eqref{eq_CMDP_Prob_Augm} are bijective using the mapping  $\eta \mapsto \Gamma(\eta)$ given in the hypothesis. Therefore, it is enough to prove that the corresponding objective and constraint costs are also equal, when respectively started in states $s$ and $(s, \1)$. 
The below theorem asserts the same  which completes the proof of Theorem \ref{thm_main_thm1} (proof is in Appendix).
\begin{thm}\label{thm_costs_eql}
For the policies $\eta\in\bar{\Pi}$ and its corresponding policy $\Gamma(\eta)\in\Pi_{\rm MR}$, the following holds:
\begin{eqnarray}
\TotCostT{\eta}{T}{i}{s}{\1} = w^{\Gamma(\eta)}_i(s) \mbox{ for all } i . \hspace{5mm} \mbox{ \eop}\label{eq_combncost_as_lincost} 
\end{eqnarray}
\end{thm}

Towards  Theorem \ref{thm_main_thm2},  observe that the problem \eqref{eq_CMDP_Prob_Augm}
 is a   standard   CMDP, but  for the  domain of optimization. LP based techniques are  well known to solve standard CMDPs  (\cite{PUT}, \cite{Eitan_CMDP}). 
 But ${\bar \Pi}$ includes only a
 %
 %
 special  class of policies which are indifferent to the augmented state $z$. This restriction is captured via the bilinear constraint of the BLP \eqref{eq_LP}, where the variables are indifferent to the value of $z$. The rest of the BLP without this particular condition is same as in \cite{Arnab_Jeffrey} extended to the constraints.  The variable $w_t(x,z,a)$ is representative of  the probability that the system is in state $(x,z)$ at time $t$ and decision $a$ is made. The rest of the proof of Theorem \ref{thm_main_thm2} is given in the Appendix. \eop

\ignore{
\begin{dmath}
{\min \sum\limits_{t=0}^{T}\sum\limits_{(x,z)\in\mathcal{\bar{X}}} \delta_t(x,z)u_t(x,z)}\\
\mbox{subject to:}\\
{u_t(x,z)\geq \bar{r}_{t,0}((x,z),a)+\mbox{\hspace{-0.2 cm}}\sum\limits_{(x^\prime,z^\prime)\mathcal{\bar{X}}}\mbox{\hspace{-0.2 cm}}\bar{Q}_t\bigg((x^\prime,z^\prime)|(x,z),a\bigg)u_{t+1}(x^\prime,z^\prime),}\\
 {\mbox{for all }t<T, (x,z)\in\mathcal{\bar{X}}, a\in\mathcal{A},}\\
 {u_T(x,z)\geq \bar{r}_{T,0}(x,z) \mbox{ for all }(x,z)\in\mathcal{\bar{X}}.}
\end{dmath}

Dual of the above LP \eqref{eq_prim_obj}-\eqref{eq_prim_cons2} along with the constraints appearing in $(P)$ is given by}
%

\section*{Conclusions}

Many applications require sequential decision models involving a combination of multiplicative and additive cost components which can be framed as CMDP.
One can convert them to standard MDP models, that was recently done (for risk-sensitive MDPs) by augmenting the state space by taking the total costs to cover the multiplicative costs, but the complexity grows exponentially in time horizon making such problems unsolvable even with moderate time horizons.  
This paper fills this gap and makes an important contribution addressing not only risk-sensitive MDPs, but also  MDPs involving a combination of multiplicative and additive cost components in objective and/or constraint. We address this by augmenting state space with binary variables such that the number of unknowns in the resulting optimization problem grow linearly in time horizon. An implementable algorithm is provided to solve such CMDPs.

\vspace{-1mm}
\section{Appendix}
We give the below two lemmas without proof. First one can be easily proved using mathematical induction and the second one is straight forward.

\begin{lemma}\label{lem_sum_prod_eql1}
Let $N\in\Z^+$ and $q\in \R^N$ be given. Then, following holds \vspace{-0.3 cm}
\begin{eqnarray}
\sum\limits_{\delta\in\{0,1\}^N}\prod\limits_{i=1}^{N} q_i^{\delta_i}(1-q_i)^{(1-\delta_i)}=1.\label{eq_general_result}\mbox{\eop}
\end{eqnarray}
\end{lemma}
\begin{lemma}\label{lem_equality_optimal}
Let $f:A\to\R$ and $g:B\to\R$ be two real valued mappings on any sets $A$ and $B$ respectively. Assume, 
i)    for each $a\in A$ there exists $b_a\in B$ such that $f(a)=g(b_a)$, and,
ii) for each $b\in B$ there exists $a_b\in A$ such that $f(a_b)=g(b)$.
Then 
$\min\limits_{a\in A} f(a) = \min\limits_{b\in B} g(b) $  and  
\begin{eqnarray*}
\mbox{if } a^*\in\argmin\limits_{a\in A} f(a) \mbox{   then }  b_{a^*}\in\argmin\limits_{b\in B} g(b), \mbox{ and similarly, } \\
\mbox{if }  b^*\in\argmin\limits_{b\in B} g(b),\mbox{   then }  a_{b^*}\in\argmin\limits_{a\in A} f(a). \hspace{16mm} \mbox{ \eop}
\end{eqnarray*}
\end{lemma}
\textbf{Proof: } Let $a^*\in\argmin\limits_{a\in A} f(a)$ and $b^{*}\in\argmin\limits_{b\in B} g(b)$, then there exists $b_{a^*}\in B$ and $a_{b^*}\in A$ such that 
\begin{eqnarray}
\min\limits_{a\in A} f(a)=f(a^*)&=&g(b_{a^*})\geq \min\limits_{b\in B} g(b)\mbox{ and }\label{eq_maxf_leq_maxg} \\
\min\limits_{a\in A} f(a) \leq f(a_{b^*})&=&g(b^*)=\min\limits_{b\in B} g(b).\label{eq_maxf_geq_maxg}
\end{eqnarray}
The inequalities \eqref{eq_maxf_leq_maxg} and \eqref{eq_maxf_geq_maxg} proves all the conclusions of the lemma.$\blacksquare$

\begin{prop}\label{prop_sum_trans_fix_ovr_z}
Fix $j\in\{0,1,\ldots,K\}$. Denote $\mathcal{Z}_j:=\{z=(z_i)_{i=0}^K\in \{0,1\}^{K+1} : z_j=1\}$. Let $z\in\mathcal{Z}_j$. 
For any $t,x,x^\prime$ and $a$, we have the following identity

\vspace{-3mm}
{\small
\begin{eqnarray}
&&\hspace{-0.6 cm}\sum\limits_{z^\prime\in\mathcal{Z}_j}\bar{Q}_t\bigg((x^\prime,z^\prime)\bigg|(x,z),a\bigg) =f_{t-1,j}(x,a,x^\prime) Q(x^\prime|x,a), \label{eq_Qt_Fix_zjzjprime_eql1}\\
&&\hspace{-1 cm}\sum\limits_{z^\prime\in\{0,1\}^{K+1}}\hspace{-0.3 cm}\bar{Q}_t\bigg((x^\prime,z^\prime)\bigg|(x,z),a\bigg) = Q(x^\prime|x,a).\label{eq_Qt_sum_zprime_eqlQ}
\end{eqnarray}

}
\end{prop}
\textbf{Proof: } Let $I:=\{i\in\{0,1,\ldots,K\}: z_i=1\}$. Clearly, $I\neq \emptyset$ because $j\in I$.
 Observe that, in the summation appearing in left hand side of \eqref{eq_Qt_Fix_zjzjprime_eql1}, the summands where $z_i=0=1-z_i^\prime$, for some $i\in\{0,1,\ldots,K\}$, is 0 and so they don't contribute to the total.  Therefore, we can restrict the summation to only those $z^\prime=(z^\prime_i)_{i=0}^K\in\mathcal{Z}_j$ which satisfy the property that $z^\prime_i=0$ for $i\notin I$. Denote $\mathcal{Z}^\prime_j :=\{z^\prime=(z^\prime_i)_{i=0}^K\in\mathcal{Z}_j:z^\prime_i=0,\mbox{ for all }i\notin I\}$. Thus, left hand side of \eqref{eq_Qt_Fix_zjzjprime_eql1} simplifies as below, where we suppress the parameters $x,a,x^\prime$ used in immediate cost functions $f_{t-1,i}$:
\ignore{
\vspace{-3mm}
{\small\begin{dmath*}
\sum\limits_{z^\prime\in\mathcal{Z}^\prime_j}\bar{Q}_t\bigg((x^\prime,z^\prime)\bigg|(x,z),a\bigg)
= {Q(x^\prime|x,a)\sum\limits_{z^\prime\in\mathcal{Z}^\prime_j}\left(\prod\limits_{i=0}^{K}(f_{t-1,i})^{z_i z_i^\prime}  (1-f_{t-1,i})^{z_i (1-z_i^\prime)}\right)}
=f_{t-1,j} Q(x^\prime|x,a)
\sum\limits_{z^\prime\in\mathcal{Z}^\prime_j} \hspace{-2mm}\left(\prod\limits_{i\in I\setminus \{j\}}\hspace{-0.2 cm}(f_{t-1,i})^{z_i z_i^\prime}(1-f_{t-1,i})^{z_i (1-z_i^\prime)}\right)
=f_{t-1,j} Q(x^\prime|x,a)
\sum\limits_{z^\prime\in\mathcal{Z}^\prime_j}\left(\prod\limits_{i\in I\setminus \{j\}}(f_{t-1,i})^{z_i^\prime}(1-f_{t-1,i})^{(1-z_i^\prime)}\right)
=f_{t-1,j} Q(x^\prime|x,a).
\end{dmath*}}
 }
 {\small
 \begin{eqnarray}
&&\sum\limits_{z^\prime\in\mathcal{Z}^\prime_j}\bar{Q}_t\bigg((x^\prime,z^\prime)\bigg|(x,z),a\bigg)\\
&=& Q(x^\prime|x,a)\sum\limits_{z^\prime\in\mathcal{Z}^\prime_j}\left(\prod\limits_{i=0}^{K}(f_{t-1,i})^{z_i z_i^\prime}  (1-f_{t-1,i})^{z_i (1-z_i^\prime)}\right)\\
&=&Q(x^\prime|x,a)\sum\limits_{z^\prime\in\mathcal{Z}^\prime_j}\left(\prod\limits_{i\in\mathcal{I}}(f_{t-1,i})^{z_i z_i^\prime}  (1-f_{t-1,i})^{z_i (1-z_i^\prime)}\right)\nonumber\\
&&\hspace{2.2 cm}\left(\prod\limits_{i\notin\mathcal{I}}(f_{t-1,i})^{z_i z_i^\prime}  (1-f_{t-1,i})^{z_i (1-z_i^\prime)}\right)\nonumber\\
&=&Q(x^\prime|x,a)\sum\limits_{z^\prime\in\mathcal{Z}^\prime_j}\left(\prod\limits_{i\in\mathcal{I}}(f_{t-1,i})^{z_i^\prime}  (1-f_{t-1,i})^{(1-z_i^\prime)}\right)\label{eq_Prod_only_I}\\
&=&Q(x^\prime|x,a)\sum\limits_{z^\prime\in\mathcal{Z}^\prime_j}(f_{t-1,j})^{z_j^\prime}  (1-f_{t-1,j})^{(1-z_j^\prime)}\nonumber\\
&&\hspace{1.9 cm}\left(\prod\limits_{i\in\mathcal{I}\setminus\{j\}}(f_{t-1,i})^{z_i^\prime}  (1-f_{t-1,i})^{(1-z_i^\prime)} \right)\\
&=&Q(x^\prime|x,a)\sum\limits_{z^\prime\in\mathcal{Z}^\prime_j}(f_{t-1,j})^{1}  (1-f_{t-1,j})^{(1-1)}\nonumber\\
&&\hspace{1.9 cm}\left(\prod\limits_{i\in\mathcal{I}\setminus\{j\}}(f_{t-1,i})^{z_i^\prime}  (1-f_{t-1,i})^{(1-z_i^\prime)} \right)\label{eq_zjprime_eqls1}\\
&=&f_{t-1,j} Q(x^\prime|x,a)
\sum\limits_{z^\prime\in\mathcal{Z}^\prime_j}\left(\prod\limits_{i\in I\setminus \{j\}}\hspace{-0.2 cm}(f_{t-1,i})^{z_i^\prime}(1-f_{t-1,i})^{(1-z_i^\prime)}\right)\nonumber\\
&=&f_{t-1,j} Q(x^\prime|x,a).
\end{eqnarray}
}

The equality in \eqref{eq_Prod_only_I} follows from its previous equation by substituting the value of $z_i$ with $z_i=1$ whenever $i\in \mathcal{I}$  and $z_i=0$ whenever $i\notin \mathcal{I}$.
The equality in \eqref{eq_zjprime_eqls1} follows from its previous equation because $\mathcal{Z}^\prime_j\subset \mathcal{Z}_j$  where   $\mathcal{Z}_j$ is defined such that $j$-th  coordinate is  1 
 and if  $z^\prime\in\mathcal{Z}^\prime_j$ then  $z^\prime\in\mathcal{Z}_j$ and thus  $z^\prime_j=1$.
 We used Lemma \ref{lem_sum_prod_eql1} in the last equality. Thus proving \eqref{eq_Qt_Fix_zjzjprime_eql1}.
 
 Following similar reasoning as used in proving \eqref{eq_Qt_Fix_zjzjprime_eql1}, we can derive the below 
{\small
 \begin{eqnarray}
     \sum\limits_{z^\prime\notin\mathcal{Z}_j}\bar{Q}_t\bigg((x^\prime,z^\prime)\bigg|(x,z),a\bigg) =(1-f_{t-1,j}(x,a,x^\prime)) Q(x^\prime|x,a)\label{eq_Qt_Fix_zjzjprime_eql0}
 \end{eqnarray}
}
Summing the equation \eqref{eq_Qt_Fix_zjzjprime_eql1} and \eqref{eq_Qt_Fix_zjzjprime_eql0} results in \eqref{eq_Qt_sum_zprime_eqlQ}.\eop
\ignore{\begin{defn}\label{def_liftpol}
Let $\pi=\{d_t\}_{t\geq 0}\in\Pi_{\rm MR}$. We call a policy $\bar{\pi}=\{\bar{d}_{t}\}_{t\geq 0} \in \bar{\Pi}_{\rm MR}$, as the `lifted' policy of $\pi$, if at every time $t=0,\ldots, T-1$,
\begin{equation}
\bar{d}_{t}(a|(x,z))=d_t(a|x), \mbox{ for all } x,z,a.
\end{equation}\eop
\end{defn}}
\begin{prop}\label{prop_Prob_diff_meas3}
 Let $\eta\in \bar{\Pi}$. For $0\leq t\leq T-1$, we have for all  $x,x^\prime,a$,
{\small
\begin{eqnarray*}
\sum\limits_{z\in\mathcal{Z}} P^{\eta}_{(s,\1)}\bigg((X_t=x,Z_t=z),A_t=a\bigg)=P^{\Gamma(\eta)}_s(X_t=x,A_t=a).
\end{eqnarray*}}
\end{prop}
\textbf{Proof: } Let $\eta=\{\bar{d}_t\}_t$ and $\Gamma(\eta)=\{d_t\}_t$. We prove the proposition using mathematical induction. At $t=0$, it is easy to prove that LHS and RHS of the above equation equals $d_0(a|s)\one_{\{x=s\}}$.
Suppose at time $t-1$, the proposition holds true. The identity at $t$ holds by following the below steps:

\vspace{-3mm}
{\small
\begin{dmath*}
 {\sum\limits_{z\in\mathcal{Z}} P^{\eta}_{(s,\1)}\bigg(X_t=x,Z_t=z,A_t=a\bigg)}={\sum\limits_{z\in\mathcal{Z}} \bar{d}_t(a|(x,z)) P^{\eta}_{(s,\1)}\bigg(X_t=x,Z_t=z\bigg)}  
 ={\sum\limits_{z\in\mathcal{Z}} d_t(a|x)P^{\eta}_{(s,\1)}\bigg((X_t=x,Z_t=z)\bigg)}
 ={d_t(a|x)\sum\limits_{x^\prime,z^\prime,a^\prime}\sum\limits_{z\in\mathcal{Z}} \bar{Q}_t\bigg((x,z)|(x^\prime,z^\prime),a^\prime\bigg)}\\{\hspace{0.5 cm}P^{\eta}_{(s,\1)}\bigg(X_{t-1}=x^\prime,Z_{t-1}=z^\prime,A_{t-1}=a^\prime\bigg)}
 ={d_t(a|x) \hspace{-0.2 cm}\sum\limits_{x^\prime,z^\prime,a^\prime}\hspace{-0.2 cm}Q(x|x^\prime,a^\prime)P^{\eta}_{(s,\1)}\bigg(X_{t-1}=x^\prime,Z_{t-1}=z^\prime,A_{t-1}=a^\prime\bigg)}
 ={d_t(a|x)\sum\limits_{x^\prime,a^\prime}Q(x|x^\prime,a^\prime) P^{\Gamma(\eta)}_s(X_{t-1}=x^\prime,A_{t-1}=a^\prime)}
 = {P^{\Gamma(\eta)}_s(X_{t}=x,A_{t}=a).}
\end{dmath*}
}
Induction hypothesis and \eqref{eq_Qt_sum_zprime_eqlQ} gives penultimate two equalities.$\blacksquare$

\ignore{
We use the following identity, whose proof is straight forward and left to the reader, to prove the proposition for $t>0$ : 
\begin{eqnarray}
\sum_{z^\prime}\bar{Q}_t((x^\prime,z^\prime)|(x,z),a) =    Q(x^\prime|x,a), \forall x^\prime,x,z,a.\label{eq_sumovrz_barQ_eqls_Q}
\end{eqnarray}
For $t>0$, simplify the summand of LHS as below,
{\small
\begin{dmath*}
{P^{\bar{\pi}}_{(s,\1)}\bigg((X_t=x,Z_t=z),A_t=a\bigg)}= {P^{\bar{\pi}}_{(s,\1)}\bigg(A_t=a\bigg|(X_t=x,Z_t=z)\bigg) P^{\bar{\pi}}_{(s,\1)}\bigg((X_t=x,Z_t=z)\bigg)}
={\bar{d}_t(a|(x,z)) P^{\bar{\pi}}_{(s,\1)}\bigg((X_t=x,Z_t=z)\bigg)}
={d_t(a|x)P^{\bar{\pi}}_{(s,\1)}\bigg((X_t=x,Z_t=z)\bigg)}
={d_t(a|x)\sum\limits_{x_i,z_i,a_i}\bigg(\bar{d}_0(a_0|(s,\1))\bar{Q}_1((x_1,z_1)|(s,\1),a_0)}\\
{\cdots \bar{d}_{t-1}(a_{t-1}|(x_{t-1},z_{t-1})) \bar{Q}_t((x,z)|(x_{t-1},z_{t-1}),a_{t-1})\bigg)}
={d_t(a|x)\sum\limits_{x_i,z_i,a_i}\bigg(d_0(a_0|s)\bar{Q}_1((x_1,z_1)|(s,\1),a_0)d_1(a_1|x)}\\
{\cdots d_{t-1}(a_{t-1}|x_{t-1}) \bar{Q}_t((x,z)|(x_{t-1},z_{t-1}),a_{t-1})\bigg).}
\end{dmath*}
}
Taking the summation over variable $z$ in above and summing backwards over the variables $z_i$ (in the product) by applying \eqref{eq_sumovrz_barQ_eqls_Q} at all times in the order $t,t-1,\ldots,1$, the LHS simplifies to
\begin{dmath*}
{\sum\limits_{z\in\mathcal{Z}} P^{\bar{\pi}}_{(s,\1)}\bigg((X_t=x,Z_t=z),A_t=a\bigg)} =
{d_t(a|x)\sum\limits_{x_i,z_i,a_i,z}\bigg(d_0(a_0|s)\bar{Q}_1((x_1,z_1)|(s,\1),a_0)d_1(a_1|x)}\\
{\cdots d_{t-1}(a_{t-1}|x_{t-1}) \bar{Q}_t((x,z)|(x_{t-1},z_{t-1}),a_{t-1})\bigg)}
={d_t(a|x)\sum\limits_{x_i,z_i,a_i}\bigg(d_0(a_0|s)\bar{Q}_1((x_1,z_1)|(s,\1),a_0)d_1(a_1|x)}\\
{\cdots d_{t-1}(a_{t-1}|x_{t-1}) Q(x|x_{t-1},a_{t-1})\bigg)}
={\sum\limits_{x_i,a_i}\bigg( d_0(a_0|s)Q(x_1|s,a_0)d_1(a_1|x)}\\
{\cdots d_{t-1}(a_{t-1}|x_{t-1}) Q(x|x_{t-1},a_{t-1})d_t(a|x)\bigg)}
={P^\pi_s(X_0=x,A_0=a).}
\end{dmath*}
}

\textbf{Proof of Theorem \ref{thm_costs_eql}:} The left hand side of \eqref{eq_combncost_as_lincost} equals,

\vspace{-4mm}
{\small 
\begin{equation}
    \TotCostT{\eta}{T}{i}{s}{\1}
=
\sum\limits_{t=0}^{T-1}\E^{\eta}_{(s,\1)}[\bar{r}_{t,i}((X_t,Z_t),A_t)]  
+ \E^{\eta}_{(s,\1)}\bigg[\bar{r}_{T,i}(X_T,Z_T)\bigg].  \label{eq_Expt_barpi_comb}
\end{equation}}
The sum of expectations in \eqref{eq_Expt_barpi_comb}  simplifies to
{\small
\begin{dmath*}
 {\sum\limits_{t=0}^{T-1}\E^{\eta}_{(s,\1)}[\bar{r}_{t,i}((X_t,Z_t),A_t)]}
= {\sum\limits_{t=0}^{T-1} \sum\limits_{x,z,a}\bar{r}_{t,i}((x,z),a) P^{\eta}_{(s,\1)}\bigg((X_t=x,Z_t=z),A_t=a \bigg) }
={\hspace{-0.1 cm}\sum\limits_{\substack{t,x,\\x^\prime,z,a }}r_{t,i}(x,a,x^\prime)Q(x^\prime|x,a)     P^{\eta}_{(s,\1)}\bigg(X_t=x,Z_t=z,A_t=a\bigg) }
={ \sum\limits_{t,x,x^\prime,a}r_{t,i}(x,a,x^\prime)Q(x^\prime|x,a)P^{\Gamma(\eta)}_s(X_t=x,A_t=a)}
= {\sum\limits_{t=0}^{T-1} \E^{\Gamma(\eta)}_s[r_{t,i}(X_t,A_t,X_{t+1})] = \E^{\Gamma(\eta)}_s\left[\sum\limits_{t=0}^{T-1} r_{t,i}(X_t,A_t,X_{t+1})\right].}
\end{dmath*}
}

\newcommand{\xv}{{\bf x}}
\newcommand{\av}{{\bf a}}
\newcommand{\zv}{{\bf z}}
\newcommand{\Q}{{\bf {\bar Q}}}
\newcommand{\bard}{{\bf {\bar D}}}
\newcommand{\dd}{{\bf {D}}}

We used Proposition \ref{prop_Prob_diff_meas3} in the third equality above. 

Define for each $i\in\{0,1,\ldots,K\}$, $\mathcal{Z}_{i}:=\{z\in\mathcal{Z}:z_i=1\}\subset \{0,1\}^{K+1}$. 
\ignore{\color{red}
Let $\xv = (x_t)_{0\leq t\le T}$
$\av = (x_t)_{0, t\le T-1}$
and $\zv = (z_t)_{0\le t\le T}$ with $z_0 = \1$, $x_0 = s$ represent the vector and we use short notation, 
$\Q_t  = {\bar Q}_t (x_t, z_t | x_{t-1}, z_{t-1}, a_{t-1})$
also let 

{\small \begin{eqnarray*}
\bard_t  (\xv,\av, \zv)  :=
\bar{d}_0(a_0|(s,\1))\Q_1 
\prod_{\tau =1}^{t}  \bar{d}_\tau (a_\tau |(x_\tau ,z_\tau ))) \Q_\tau
\end{eqnarray*}
With these notations we have:

{\small \begin{eqnarray*}\E^{\eta}_{(s,\1)}\bigg[\bar{r}_{T,i}(X_T,Z_T)\bigg] &=&\alpha_i \E^{\eta}_{(s,\1)}[Z_{T,i}]={\alpha_i P^{\eta}_{(s,\1)}[Z_{T,i}=1]}    \\
&=& \alpha_i \sum_{\xv,\zv, \av} \bard_{T-1}  (\xv,\av, \zv)
\end{eqnarray*}}
}
}
We denote the decisions of the policy $\eta$ by $\{\bar{d}_t\}_t$ and that of policy $\Gamma(\eta)$ by $\{d_t\}_t$. To simplify the terminal cost in \eqref{eq_Expt_barpi_comb}, we sum the sample path probabilities:

{\small
\begin{dmath*}
{\E^{\eta}_{(s,\1)}\bigg[\bar{r}_{T,i}(X_T,Z_T)\bigg]=\alpha_i \E^{\eta}_{(s,\1)}[Z_{T,i}]=\alpha_i P^{\eta}_{(s,\1)}[Z_{T,i}=1]}
= {\alpha_i\hspace{-0.2 cm}\sum\limits_{\substack{x_t\in\mathcal{X},\\z_t\in\mathcal{Z}_i,\\a_t\in\mathcal{A}}}\hspace{-0.1 cm}\prod\limits_{t=0}^{T-1} \bar{d}_t(a_t|(x_t,z_t))\bar{Q}_{t+1}((x_{t+1},z_{t+1})|(x_{t},z_{t}),a_{t})}
={\alpha_i\hspace{-0.2 cm}\sum\limits_{\substack{x_t\in\mathcal{X},\\z_t\in\mathcal{Z}_i,\\a_t\in\mathcal{A}}}\hspace{-0.1 cm}\prod\limits_{t=0}^{T-1} d_t(a_t|x_t)\bar{Q}_{t+1}((x_{t+1},z_{t+1})|(x_{t},z_{t}),a_{t})}
={\alpha_i\hspace{-0.2 cm}\sum\limits_{\substack{x_t\in\mathcal{X},\\a_t\in\mathcal{A}}}d_t(a_t|x_t)Q(x_{t+1}|x_{t},a_{t})f_{t,i}(x_{t},a_{t},x_{t+1})}
=\alpha_i \E^{\Gamma[\eta]}_s\left[\prod\limits_{t=0}^{T-1} f_{t,i}(X_{t},A_{t},X_{t+1})\right].\label{eq_multplcost_secondterm}
\end{dmath*}}

\ignore{
\begin{dmath*}
\E^{\eta}_{(s,\1)}\bigg[\bar{r}_{T,i}(X_T,Z_T)\bigg]=\alpha_i \E^{\eta}_{(s,\1)}[Z_{T,i}]={\alpha_i P^{\eta}_{(s,\1)}[Z_{T,i}=1]}
= \alpha_i\hspace{-0.2 cm}\sum\limits_{\substack{x_t\in\mathcal{X},\\z_t\in\mathcal{Z}_i,\\a_t\in\mathcal{A}}}\hspace{-0.1 cm} \bar{d}_0(a_0|(s,\1))\bar{Q}_1((x_1,z_1)|(s,\1),a_0)\bar{d}_1(a_1|(x_1,z_1)))\\
\mbox{ \hspace{0.7 cm}  }\bar{Q}_2((x_2,z_2)|(x_1,z_1),a_1)
\cdots \bar{d}_{T-1}(a_{T-1}|(x_{T-1},z_{T-1}))\\
\mbox{ \hspace{0.7 cm}  }\bar{Q}_{T}((x_T,z_T)|(x_{T-1},z_{T-1}),a_{T-1})
=\alpha_i\hspace{-0.2 cm}\sum\limits_{\substack{x_t\in\mathcal{X},\\z_t\in\mathcal{Z}_i,\\a_t\in\mathcal{A}}}\hspace{-0.1 cm}d_0(a_0|s) \bar{Q}_1((x_1,z_1)|(s,\1),a_0)d_1(a_1|x_1)\\ 
\mbox{ \hspace{0.7 cm}  }\bar{Q}_2((x_2,z_2)|(x_1,z_1),a_1)\cdots d_{T-1}(a_{T-1}|x_{T-1})\\
\mbox{ \hspace{0.7 cm}  }\bar{Q}_{T}((x_T,z_T)|(x_{T-1},z_{T-1}),a_{T-1})
= \alpha_i\hspace{-0.2 cm}\sum\limits_{\substack{x_t\in\mathcal{X},\\z_t\in\mathcal{Z}_i,\\a_t\in\mathcal{A}}}\hspace{-0.1 cm} d_0(a_0|s)\bar{Q}_1((x_1,z_1)|(s,\1),a_0)d_1(a_1|x_1)\\
\mbox{ \hspace{0.7 cm}  }\bar{Q}_2((x_2,z_2)|(x_1,z_1),a_1) \cdots d_{T-1}(a_{T-1}|x_{T-1})\\
\mbox{ \hspace{0.7 cm}  }{\color{blue}\bigg(}\sum\limits_{z_T\in\mathcal{Z}_i, z_{T,i}=1}\bar{Q}_{T}((x_T,z_T)|(x_{T-1},z_{T-1}),a_{T-1}){\color{blue}\bigg)}
= \alpha_i\hspace{-0.2 cm}  \sum\limits_{\substack{x_t\in\mathcal{X},\\z_t\in\mathcal{Z}_i,\\a_t\in\mathcal{A}}}\hspace{-0.1 cm} d_0(a_0|s)\bar{Q}_1((x_1,z_1)|(s,\1),a_0)d_1(a_1|x_1)\\
\mbox{ \hspace{0.7 cm}  }\bar{Q}_2((x_2,z_2)|(x_1,z_1),a_1)\cdots d_{T-1}(a_{T-1}|x_{T-1}) \\
\mbox{ \hspace{0.7 cm}  }Q(x_{T}|x_{T-1},a_{T-1}) f_{T-1,i}(x_{T-1},a_{T-1},x_{T})
=\alpha_i\hspace{-0.2 cm}  \sum\limits_{\substack{x_t\in\mathcal{X},\\a_t\in\mathcal{A}}}\hspace{-0.1 cm} d_0(a_0|s) Q(x_{1}|s,a_{0}) f_{0,i}(s,a_{0},x_{1}) d_1(a_1|x_1) \\
\mbox{ \hspace{0.7 cm} } Q(x_{2}|x_1,a_{1})f_{1,i}(x_{1},a_{1},x_{2})\cdots d_{T-1}(a_{T-1}|x_{T-1})\\
\mbox{ \hspace{0.7 cm} }Q(x_{T}|x_{T-1},a_{T-1})f_{T-1,i}(x_{T-1},a_{T-1},x_{T})
=\alpha_i \E^\pi_s\left[\prod\limits_{t=0}^{T-1} f_{t,i}(X_{t},A_{t},X_{t+1})\right].\label{eq_multplcost_secondterm}
\end{dmath*}
}
The penultimate equality above is due to successive application of \eqref{eq_Qt_Fix_zjzjprime_eql1} backwards in time for $t=T,T-1,\ldots,1$.
Replacing back each of the expectation operation terms in \eqref{eq_Expt_barpi_comb} proves the identity \eqref{eq_combncost_as_lincost}.\eop

\ignore{\textbf{Proof of Theorem \ref{thm_main_thm2}: } Denote the set of all MR policies in the augmented model $M$ by $\bar{\Pi}_{\rm MR}$. Recall, $\bar{\Pi}\subset \bar{\Pi}_{\rm MR}$ is a set of MR policies that are indifferent to the $z$ coordinate. The problem $(P^\prime):\min_{\eta\in\bar{\Pi}_{\rm MR}}\TotCostT{\eta}{T}{i}{s}{\1}$ is a standard linear MDP. The LP formulation given in \eqref{eq_LP} without the last two constraints gives the solution to the problem $(P^\prime)$ (\cite{Arnab_Jeffrey}). The penultimate constraint of the LP corresponds to the constraints of the form $\TotCostT{\eta}{T}{i}{s}{\1}\leq b_i$ for all $1\leq i\leq K$ along with the cost constraints for the standard CMDP problem. 
The problem \eqref{eq_CMDP_Prob_Augm} has a restricted domain $\bar{\Pi}$ which lies in $ \bar{\Pi}_{\rm MR}$. This restriction is captured via the last constraint of the LP \eqref{eq_LP} {\color{blue} (relationship between feasible point in LP and policy in $\bar{\Pi}$ has to be established here}). Thus, if $v_t^*(x,z,a)$ is a solution of the LP, then the optimal policy $\eta^*=\{d_t^*\}_t$ for the problem \eqref{eq_CMDP_Prob_Augm} is given by,
\begin{eqnarray*}
d_t^*(a|(x,z)):=\frac{v_t^*(x,z,a)}{\sum\limits_{a^\prime}v_t^*(x,z,a^\prime)}=\frac{v_t^*(x,\1,a)}{\sum\limits_{a^\prime}v_t^*(x,\1,a^\prime)}.
\end{eqnarray*}
Thus the policy, $\pi(\eta^*)=\pi^*$ given in \eqref{eq_sol_orig_prb}  is the optimal policy for the problem \eqref{eq_CMDP_Prob}.}

\textbf{Proof of Theorem \ref{thm_main_thm2}: }

Let $\bar{\Pi}_{\rm MR}$ denote space of all MR policies in the MDP $\bar{M}$. Recall $\bar{\Pi}$ is a set of policies that are indifferent to augmented component $z$. Note that $\bar{\Pi}\subset \bar{\Pi}_{\rm MR}$. Denote the feasible region of the BLP by $\mathcal{Q}$. 
Observe that objective function and all the constraints except the bilinear constraints are linear and therefore the problem BLP without the bilinear constraints is indeed a LP that solves the CMDP \eqref{eq_CMDP_Prob_Augm} with the domain $\bar{\Pi}_{\rm MR}$, instead of $\bar{\Pi}$. Denote the feasible region of this LP by $\mathcal{L}$. Clearly $\mathcal{Q}\subseteq\mathcal{L}$.

We first claim that the feasible region of the problem \eqref{eq_CMDP_Prob_Augm} is bijective to  $\mathcal{Q}$	
by defining the mappings  $\pi\mapsto w_\pi$ and $w\mapsto \pi_w$ respectively between these two sets.

It is easy to see that, given a feasible vector  $w=\{\Wv_1, \Wv_2\} \in\mathcal{Q}$, constructing the policy $\pi_w=\{d_t\}_t$ by rationalising over $a$ and applying the bilinear constraints immediately as below
\begin{eqnarray*}
d_t(a|(x,z))&:=&\frac{w_t(x,z,a)}{\sum\limits_{a^\prime}w_t(x,z,a^\prime)}=\frac{w_t(x,\1,a)}{\sum\limits_{a^\prime}w_t(x,\1,a^\prime)}\\
&=&d_t(a|(x,\1)),
\end{eqnarray*}
makes $\pi_w\in\bar{\Pi}$.
Also we know that, for a given feasible policy of the (augmented) problem \eqref{eq_CMDP_Prob_Augm}, say, $\pi=\{d_t\}_t\in\bar{\Pi}\subset  \bar{\Pi}_{\rm MR}$, there exists a feasible vector  $w_\pi = \{\Wv_1, \Wv_2\} \in\mathcal{L}$ with $w_t(x,z,a):=P^\pi_s(X_t=x,Z_t=z,A_t=a)$ (see \cite{PUT}). Now, to prove $w_\pi\in\mathcal{Q}$,  it is enough to prove that $w_\pi$ satisfies the bilinear constraints. 

From the literature on linear MDPs   (\cite{PUT,Kavitha_NH_Atul, Arnab_Jeffrey}), we know that, 
the mappings $w\mapsto \pi_w$ and $\pi\mapsto w_\pi$ are such that $w_{\pi_w}=w$ and $\pi_{w_\pi}=\pi$.

Choose any arbitrary $z\in\mathcal{Z}$, a $K+1$ dimensional vector of 0s and 1s. Then, the $t$-th decision w.r.t the policy $\pi_{w_\pi}$ for the two states $(x,z_1),(x,\1)\in\mathcal{\bar{X}}$ is given by

\vspace{-3mm}
{\small\begin{equation}
d_t(a|(x,z))= \frac{(w_\pi)_t(x,z,a)}{\sum\limits_{a^\prime}(w_\pi)_t(x,z,a^\prime)},\ 
d_t(a|(x,\1)) = \frac{(w_\pi)_t(x,\1,a)}{\sum\limits_{a^\prime}(w_\pi)_t(x,\1,a^\prime)}.
\end{equation}}
Because $\pi_{w_\pi}=\pi\in\bar{\Pi}$, the decisions are indifferent to the vector $z$ and $\1$, that is, $d_t(a|(x,z))=d_t(a|(x,\1))$ implying, 
$(w_{\pi})_{t}(x, z, a) \sum_{a'} (w_{\pi})_{t} (x, \1, a') = (w_{\pi})_{t} (x, \1, a) \sum_{a'} (w_{\pi})_{t} (x, z, a')$. Thus satisfying the bilinear constraint.

The below are two claims that suffices to complete the proof.

\textbf{Claim 1:} For a given feasible policy $\pi$ of the augmented problem \eqref{eq_CMDP_Prob_Augm}, the objective function of \eqref{eq_CMDP_Prob_Augm} when evaluated at $\pi$ is equal the objective function of the BLP (almost LP except for bilinear constraint) when evaluated at its feasible point $w_\pi$.

\textbf{Claim 2: } For a given feasible vector $w=(\Wv_1,\Wv_2)$ of the BLP,  the objective function of the BLP evaluated at $w$ is equal to the objective function of augmented problem \eqref{eq_CMDP_Prob_Augm} when evaluated at $\pi_w$.

The above two claims complete the proof of the theorem, by applying the Lemma \ref{lem_equality_optimal}.
\eop

\textbf{\textit{Proof of the claims 1 and 2 in the above proof.}}

Let $\pi\in \bar{\Pi}\subset  \bar{\Pi}_{\rm MR}$ be arbitrary. The corresponding feasible point of the BLP is $w_\pi=\{\Wv_1, \Wv_2\} $ defined by $w_t(x,z,a):=P^{\pi}_{(s,\1)}((X_t=x,Z_t=z),A_t=a)$ for $t<T$ and $w_T(x,z):=P^{\pi}_{(s,\1)}((X_T=x,Z_T=z))$. The objective of \eqref{eq_CMDP_Prob_Augm} under $\pi$ is
{\small\begin{eqnarray}
\TotCostT{\pi}{T}{0}{s}{\1} &:=& \E^{\pi}_{(s,\1)}\bigg[\sum\limits_{t=0}^{T-1}\bar{r}_{t,0}\bigg((X_t,Z_t),A_t\bigg)  + \bar{r}_{T,0}(X_T,Z_T)\bigg]\nonumber\\
&=&  \E^{\pi}_{(s,\1)}\bigg[\sum\limits_{t=0}^{T-1}\bar{r}_{t,0}\bigg((X_t,Z_t),A_t\bigg)\bigg] \nonumber\\
&& + \E^{\pi}_{(s,\1)}\bigg[\bar{r}_{T,0}(X_T,Z_T)\bigg].\label{eq_aug_cost_policy_splitexpt}
\end{eqnarray}}
Simplifying the first expectation in \eqref{eq_aug_cost_policy_splitexpt}, as below
{\small 
\begin{dmath*}
{\E^{\pi}_{(s,\1)}\bigg[\sum\limits_{t=0}^{T-1}\bar{r}_{t,0}\bigg((X_t,Z_t),A_t\bigg)\bigg] =
 \sum\limits_{t=0}^{T-1}\E^{\pi}_{(s,\1)}[\bar{r}_{t,0}((X_t,Z_t),A_t)]}
= {\sum\limits_{t=0}^{T-1} \sum\limits_{x,z,a}\bar{r}_{t,0}((x,z),a) P^{\pi}_{(s,\1)}\bigg((X_t=x,Z_t=z),A_t=a \bigg) }
={\sum\limits_{t=0}^{T-1} \sum\limits_{x,z,a}\bar{r}_{t,0}((x,z),a) w_t(x,z,a).}
\end{dmath*}}
Simplifying the second expectation in \eqref{eq_aug_cost_policy_splitexpt}, as below
{\small\begin{dmath*}
{\E^{\pi}_{(s,\1)}\bigg[\bar{r}_{T,0}(X_T,Z_T)\bigg] = \sum\limits_{x,z}\bar{r}_{T,0}((x,z))P^{\pi}_{(s,\1)}((X_T=x,Z_T=z)) }
={\sum\limits_{x,z}\bar{r}_{T,0}((x,z))w_T(x,z).}
\end{dmath*}}
Summing the two expectations gives the linear objective of BLP problem for the variable vector $w_\pi$ for an arbitrary $\pi$. Thus proving the point no 1.

To prove claim 2, let $w=\{\Wv_1, \Wv_2\}$ be a feasible point in the BLP. Let $\pi_w$ be the corresponding feasible point to the problem \eqref{eq_CMDP_Prob_Augm}. We know that $w_{\pi_w}=w$ and claim 1, objective function of \eqref{eq_CMDP_Prob_Augm} evaluated at policy $\pi_w$ is equal to the objective function of BLP evaluated at $w_{\pi_w}=w$. This proves claim 2.\eop

\addtolength{\textheight}{-12cm}   





\end{document}